\documentclass[10pt]{article}

\usepackage{geometry}                		
\usepackage{graphicx}	
\usepackage{xcolor}				
\usepackage{amssymb}
\usepackage{amsmath}
\usepackage{amsfonts}
\usepackage{mathtools}
\usepackage[utf8]{inputenc}
\usepackage{amsthm}

\usepackage{amsthm}
\newtheorem{theorem}{Theorem}

\usepackage{algorithmic}
\usepackage{algorithm}
\usepackage{booktabs}
\usepackage{ulem}
\usepackage{authblk}

\title{Spectral and norm estimates for matrix sequences arising from a finite difference approximation of elliptic operators}

\date{}

\author[1]{Armando Coco}
\author[2]{Sven-Erik Ekström}
\author[3]{Giovanni Russo}
\author[2,4]{Stefano Serra-Capizzano}
\author[3,5]{Santina Chiara Stissi}
\affil[1]{\small School of Engineering, Computing and Mathematics, Oxford Brookes University - OX33 1HX, Oxford (UK)}
\affil[2]{\small Department of Information Technology, Division of Scientific Computing, Uppsala University - ITC, Lägerhyddsv. 2, hus 2,  P.O. Box 337, SE-751 05, Uppsala (SWEDEN)}
\affil[3]{\small Department of Mathematics and  Informatics,  Catania  University - Viale A. Doria 6, 95125 Catania  (ITALY)}
\affil[4]{\small Department of Humanities and Innovation, Insubria University - via Valleggio 11,  22100 Como (ITALY)}
\affil[5]{\small Istituto Nazionale di Geofisica e Vulcanologia - Piazza Roma 2,  95125 Catania (ITALY)}

\begin{document}

\maketitle

\begin{abstract}
When approximating elliptic problems by using specialized approximation techniques, we obtain large structured matrices whose analysis provides information on the stability of the method. Here we provide spectral and norm estimates for matrix sequences arising from the approximation of the Laplacian via ad hoc finite differences. The analysis involves several tools from matrix theory and in particular from the setting of Toeplitz operators and Generalized Locally Toeplitz matrix sequences. Several numerical experiments are conducted, which confirm the correctness of the theoretical findings.
\end{abstract}

{\bf Keywords:} Toeplitz matrix, generating function and spectral symbol, approximation of differential operators.

\section{Introduction}\label{sec:intro}

In the numerical approximation of elliptic differential equations, by using specialized approximation techniques, we obtain large structured matrices whose analysis provides information on the stability of the method. Here we provide spectral and norm estimates for matrix sequences arising from the approximation of the Laplacian via ad hoc finite differences that is from the Coco--Russo method~\cite{coco2013finite}.

The analysis involves several tools from matrix theory and in particular from the setting of Toeplitz operators and Generalized Locally Toeplitz (GLT) matrix sequences. Several numerical experiments are conducted, which confirm the theoretical findings.

The paper is organized as follows. Subsection \ref{sec:intro-method} contains a motivation and a description of the Coco--Russo method, together with a brief account on the related literature. Subsection \ref{sec:intro-toeplitz} contains the necessary tools from the Toeplitz technology, while Section \ref{sec:pb 1d} contains the matrix formulation in $1\mathbb{D}$ in the language of Toeplitz structures, the analysis of the norm estimates in $1\mathbb{D}$, together with related numerical experiments and a preliminary discussion on the spectral features of the involved matrix-sequences. Section \ref{sec:2D} contains more details on the $2\mathbb{D}$ method, on its matrix formulation, on the spectral results in $1\mathbb{D}$ and in $2\mathbb{D}$, and the basic tools taken the GLT theory. A discussion on the more challenging case of the norm estimates in $2\mathbb{D}$ is also provided. A conclusion section ends the paper with a mention to a few open problems.

\subsection{Method description and motivation}\label{sec:intro-method}

The design of numerical methods to solve Partial Differential Equations (PDE) on complex-shaped domains is obtaining an increasing interest in the scientific community. One of the bottlenecks of modern computer simulations is the modelling of physical processes around $3\mathbb{D}$ complex-shaped objects through PDE.
Finite Element Methods (FEM) are well-established approaches to solve PDE and supported by rigorous theoretical analysis developed in the last decades to prove the convergence and accuracy order of the method when the grid size approaches zero.

However, some critical limitations are commonly associated in literature with FEM, especially when applied to curved boundaries.
 In particular, the generation of elements to conform highly varying curvatures of the boundary might become cumbersome, especially if the domain changes its shape over time. Also, the design of a balanced partition of the mesh for parallel FEM is unhandy. For these reasons, approaches based on Finite Difference Methods (FDM) where the domain is immersed into a fixed grid are increasing their popularity in literature, since they do not require any mesh generation effort and at the same time allow for a natural design of parallel solvers.

On the other hand, FDM are commonly based on heuristic approaches and convergence and stability analysis are not sufficiently developed in literature, especially for the case of curved boundaries.

The Immersed Boundary Method proposed by Peskin in~\cite{peskin1977numerical} and further developed by LeVeque and Li in~\cite{leveque1994immersed} is a pioneer approach based on FDM for general domains immersed on fixed grids.

A more recent approach is the Ghost-Fluid Method proposed by Fedkiw et al.~in~\cite{fedkiw1999non} and further extended to higher accuracy by Gibou et al.~in~\cite{gibou2002second, gibou2005fourth}, where the values on grid nodes just outside the domain (ghost points) are obtained by accurate extrapolations of the boundary condition from inside values.

In~\cite{coco2013finite}, the authors present a highly efficient and accurate ghost-point method to solve a Poisson equation on a complex-shaped domain, modelled by a level-set function. Several numerical tests were presented to confirm the accuracy order and the efficiency of the multigrid solver. However, a theoretical analysis was missing. The method has been extended to several applications, such as compressible fluids in moving domains~\cite{chertock2018second} or volcanology~\cite{coco2014second}.

In this paper we present a technique to prove the stability of the Coco--Russo method~\cite{coco2013finite} and the convergence to the predicted order of accuracy.

\begin{figure}[htbp]
\centering
\includegraphics[scale=.85]{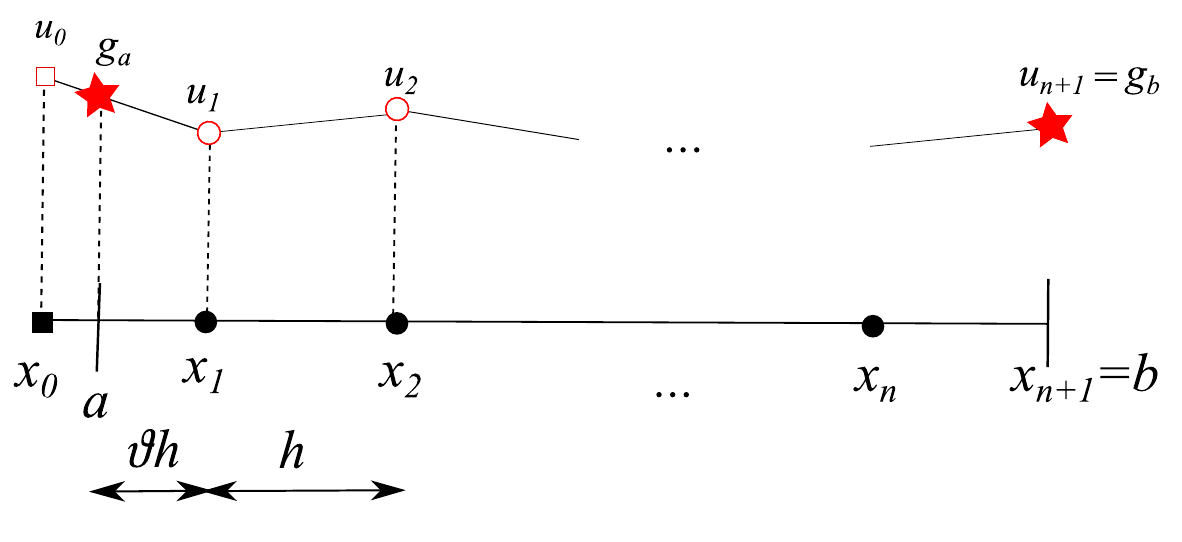}
\caption{Discretization of the $1\mathbb{D}$ domain. Full black circles are the interior grid points, while the full square is the ghost point. Boundary values are indicated with red stars. Linear extrapolation is used to define the ghost value $u_0$ from $u_1$ and the left boundary value $g_a$. }
\label{domain1D}
\end{figure}

We start from the $1\mathbb{D}$ problem.
Consider the elliptic boundary-value problem:
\begin{eqnarray}\label{mainProblem}
- \Delta u = f \text{ on } \Omega = (a,b) \\
u(a) = g_a, \quad u(b)=g_b
\end{eqnarray}
and a one-dimensional uniform grid $\mathcal{G}_h = \left\{ x_0, x_1, \ldots, x_{n+1} \right\}$ with a constant spatial step $h = x_{i}-x_{i-1}$, for $i=1,\ldots,n+1$.
Then, $x_i = x_0 +  i\, h$.
Let $x_0<a<x_1$ and $x_{n+1} = b$ (see Fig.~\ref{domain1D}).

The elliptic equation $- \Delta u = f$ is discretized by central differences on $x_i$ for $i=1,\ldots, n$ and the boundary condition on $x=b=x_{n+1}$ is included in the internal discretization (this is the so called {\em eliminated boundary condition} approach):
\[
\frac{-u_{i-1}+2u_i- u_{i+1}}{h^2} = f_i \text{ for } i=1,\ldots,n-1,
\]
\[
\frac{2u_{n} - u_{n-1}}{h^2} = f_{n}+\frac{g_{b}}{h^2}.
\]
The boundary condition on $x=a$ is approximated by $q(a) = g_a$, where $q(x)$ is the polynomial of degree $s-1$ that interpolates $u$ on the grid points $u_0, u_{1}, \ldots, u_{s-1}$. We call $s$ the stencil size for the boundary condition on $x=a$.

The discretization of the boundary condition can be represented as:
\begin{equation}\label{BC1D}
\sum_{i=0}^{s-1} c_{i} u_{i} = g_\alpha.
\end{equation}
For $s=2$ we have
\begin{equation}\label{BC1D_Dir}
\vartheta u_0 + (1-\vartheta) u_1 = g_\alpha,
\end{equation}
where $\vartheta = (x_1-a)/h$.
The grid point $x_0$ is called ghost point and $u_0$ is the ghost value.

Although we can follow a similar technique for the boundary condition on $x=a$ to the one that we adopted for $x=b$ (i.e.~we can solve~\eqref{BC1D_Dir} for $u_0$ and substitute its value into the internal equation for $x_1$),
we keep a non-eliminated boundary condition approach in order to develop a theoretical analysis that can be straightforwardly extended to higher dimensional cases, where the eliminated approach is impractical.

The discretized problem is then a linear system $A_h \textbf{u}_h = \textbf{f}_h$ where $A_h \in \mathbb{R}^{(n+1) \times (n+1)}$:
\begin{align}\label{LS1D}
A_h\mathbf{u}_h=
\begin{bmatrix*}[c]
\vartheta&1-\vartheta&\phantom{-}0&\ldots&\ldots&\phantom{-}0\\
-\frac{1}{h^2}&\frac{2}{h^2}&-\frac{1}{h^2}&\phantom{-}0&\ldots&\phantom{-}0\\
\phantom{-}0&\ddots&\ddots&\ddots&\ddots&\vdots\\
\vdots&\ddots&\ddots&\ddots&\ddots&\phantom{-}0\\
\phantom{-}0&\ldots&\phantom{-}0&-\frac{1}{h^2}&\phantom{-}\frac{2}{h^2}&-\frac{1}{h^2}\\
\phantom{-}0&0&\ldots&\phantom{-}0&-\frac{1}{h^2}&\phantom{-}\frac{2}{h^2}
\end{bmatrix*}
\begin{bmatrix*}
u_0\\u_1\\\vdots\\\vdots\\u_{n-1}\\u_{n}
\end{bmatrix*}
=
\begin{bmatrix*}
g_a\\f_1\\f_2\\\vdots\\f_{n-1}\\f_{n}+\frac{g_b}{h^2}
\end{bmatrix*}=\mathbf{f}_h,
\end{align}
where $h=(b-x_0)/(n+1)$.

\subsection{Toeplitz structures and related tools}\label{sec:intro-toeplitz}
Let $T_n$ be a Toeplitz matrix of order $n$ and let $\omega<n$ be a positive integer
\begin{align}
\setlength{\arraycolsep}{0.15em}
T_n={\footnotesize\left[\begin{array}{cccccc}
a_0 & \cdots & a_{-\omega}&  &  \\[-0.5em]
\vdots & \ddots   &  & \ddots &  \\[-0.5em]
a_{\omega} & & \ddots &     & \ddots \\[-0.5em]
 &\ddots &   &\ddots   &  & a_{-\omega}\\[-0.5em]
 & & \ddots &  & \ddots & \vdots \\[-0.5em]
 & &  & a_{\omega} &\cdots & a_0
\end{array}\right],}
\label{toeplitz-m}
\end{align}
where the coefficients $a_k$, $k=-\omega,\ldots,\omega$, are complex numbers.

Let $f\in L^1(-\pi,\pi)$ and let $T_n(f)$ be the Toeplitz matrix generated by $f$ i.e. $\left(T_n(f)\right)_{s,t}=a_{s-t}(f)$, $s,t=1,\ldots,n$, with $f$ indicated as \emph{generating function}  of $\{T_n(f)\}$  and with $a_k(f)$ being the $k$-th Fourier coefficient of $f$ that is

\begin{equation}\label{fou}
a_k(f)=\frac{1}{2\pi}\int_{-\pi}^{\pi} f(\theta)\mathrm{e}^{-\mathbf{i}k\theta}\, \mathrm{d}\theta,\quad \mathbf{i}^2=-1,\ k\in \mathbb{Z}.
\end{equation}
With these notations the matrix reported in (\ref{toeplitz-m}) can be written as $T_n=T_n(f)$, where the generating function is 
$f(\theta)=\sum_{j=-\omega}^{\omega} a_j \mathrm{e}^{\mathbf{i}j\theta}$. It is worth noticing that study of the generating function gives plenty of information on the spectrum of $T_n(f)$ for any fixed $n$, and also asymptotically as the matrix-size $n$ diverges to infinity (see \cite{MR3674485} and \cite{GLTvol2} for the multilevel setting). For instance,r if $f$ is real-valued almost everywhere (a.e.), then $T_n(f)$ is Hermitian for all $n$. Furthermore, when $f$ is real-valued and even a.e., the matrix $T_n(f)$ is (real) symmetric for all $n$, while $f$ real-valued and nonnegative a.e., but not identically zero a.e., implies that $T_n(f)$ is Hermitian positive definite for all $n$: in such a setting the considered matrix-sequence could be ill-conditioned and indeed if $f$ is nonnegative and bounded with essential supremum equal to $M>0$ and a unique zero of order $\alpha>0$, then the maximal eigenvalue converges monotonically from below to $M$, whereas the minimal eigenvalues converges to zero monotonically from above with a speed dictated by $\alpha$, that is the minimal eigenvalue is asymptotical to $n^{-\alpha}$.
In many practical applications we remind that it is required to solve numerically linear systems of Toeplitz kind and of (very) large dimensions and hence several specialized techniques of iterative type, such as preconditioned Krylov methods and ad hoc multigrid procedures have been designed; we refer the interested reader to the books \cite{MR2108963,MR2376196} and to the references therein. We recall that such types of large Toeplitz linear systems  emerge from specific  applications involving e.g. the numerical solution of (integro-) differential equations and of problems with Markov chains.

\section{Matrix formulation and notation in $1\mathbb{D}$}\label{sec:pb 1d}
\noindent The linear system to solve is \eqref{LS1D},
and we can decompose the matrix $A_h\in\mathbb{R}^{(n+1)\times(n+1)}$ as follows
\begin{align}
A_h&=\frac{1}{h^2}
\begin{bmatrix*}
\phantom{-}2&-1&\phantom{-}0&\phantom{-}\ldots&\phantom{-}\ldots&\phantom{-}0\\
-1&\phantom{-}2&-1&\phantom{-}0&\phantom{-}\ldots&\phantom{-}0\\
\phantom{-}0&\phantom{-}\ddots&\phantom{-}\ddots&\phantom{-}\ddots&\phantom{-}\ddots&\phantom{-}\vdots\\
\phantom{-}\vdots&\phantom{-}\ddots&\phantom{-}\ddots&\phantom{-}\ddots&\phantom{-}\ddots&\phantom{-}0\\
\phantom{-}0&\phantom{-}\ldots&\phantom{-}0&-1&\phantom{-}2&-1\\
\phantom{-}0&\phantom{-}0&\phantom{-}\ldots&\phantom{-}0&-1&\phantom{-}2
\end{bmatrix*}
+
\frac{1}{h^2}\begin{bmatrix*}
\vartheta h^2-2&(1-\vartheta)h^2+1&0&\ldots &0\\
0&0&0&\ldots&0\\
\vdots&\vdots&\vdots&&\vdots\\
\vdots&\vdots&\vdots&&\vdots\\
\vdots&\vdots&\vdots&&\vdots\\
0&0&0&\ldots&0\\
\end{bmatrix*}\nonumber \\
&=\frac{1}{h^2}T_{n+1}(2-2\cos(\theta))+\frac{1}{h^2}
\begin{bmatrix*}
\mathbf{v}_h^{\mathrm{T}}\label{eq:Tn}\\
\mathbf{0}\\
\vdots\\
\mathbf{0}
\end{bmatrix*}\\
&=S_{n+1}+\frac{1}{h^2}\mathbf{e}_1 \mathbf{v}_h^{\mathrm{T}},\label{eq:Sn}
\end{align}
where $T_{n+1}(f)$ in \eqref{eq:Tn} is the Toeplitz matrix generated by $f$ according to (\ref{fou}), with $f(\theta)=2-2\cos(\theta)$ so that,
in the matrix in  (\ref{toeplitz-m}), we have $\alpha=1$, $a_0=2,a_1=a_{-1}=-1$. Furthermore we have defined $S_{n+1}$ in \eqref{eq:Sn} as
\[
S_{n+1}=\frac{1}{h^2}T_{n+1}(2-2\cos(\theta)).
\]
For this matrix everything is known and in fact
\[
T_{n+1}(2-2\cos(\theta))=QDQ
\]
with $Q$ real symmetric and orthogonal and
\[
Q=Q_{n+1}=\left(\sqrt{\frac{2}{n+2}}\sin\left(\frac{st\pi}{n+2}\right)\right)_{s,t=1}^{n+1},\ \ \
D={\rm diag}\left(4\sin^2\left(\frac{s\pi}{2(n+2)}\right)\right).
\]
Hence its conditioning $\kappa_2(\cdot)$ in spectral norm (the one induced by the Euclidean vector norm) is exactly known and it is equal to
\[
\kappa_2(S_{n+1})=\sin^2\left(\frac{(n+1)\pi}{2(n+2)}\right)\sin^{-2}\left(\frac{\pi}{2(n+2)}\right)\approx
\frac{4}{\pi^2} n^2,
\]
where $a_n\approx b_n$ means $a_n=b_n(1+o(1))$ and where, in our setting,  a even more precise relation can be derived, that is that is $\kappa_2(S_{n+1})=\frac{4}{\pi^2} n^2 +O(1)$. Since everything is known regarding the term $S_{n+1}$ our idea is to reduce the analysis as much as possible to information concerning the matrix $S_{n+1}$ and its inverse and to this end the application of the Sherman--Morrison--Woodbury is appropriate.

The Sherman--Morrison--Woodbury formula states that for and invertible square matrix $A$, column vectors $\mathbf{u}$ and $\mathbf{v}$, and $1+\mathbf{v}^{\mathrm{T}}A^{-1}\mathbf{u}\neq 0$
\begin{align}
\left(A+\mathbf{u}\mathbf{v}^{\mathrm{T}}\right)^{-1}=A^{-1}-\frac{A^{-1}\mathbf{u}\mathbf{v}^{\mathrm{T}}A^{-1}}{1+\mathbf{v}^{\mathrm{T}}A^{-1}\mathbf{u}}.
\label{eq:smw}
\end{align}
and thus we can obtain in our setting defined above in \eqref{eq:smw} with $A=S_{n+1}$ and $\mathbf{u}=\frac{\mathbf{e_1}}{h^2}$ and $\mathbf{v}=\mathbf{v}_h$.
\[
\left(S_{n+1}+\frac{1}{h^2}\mathbf{e}_1 \mathbf{v}_h^{\mathrm{T}}\right)^{-1}=S_{n+1}^{-1}-\frac{S_{n+1}^{-1}\frac{1}{h^2}\mathbf{e}_1 \mathbf{v}_h^{\mathrm{T}}S_{n+1}^{-1}}{1+\mathbf{v}_h^{\mathrm{T}}S_{n+1}^{-1}\frac{\mathbf{e}_1}{h^2}}.
\]
or
\begin{align}
A_h^{-1}&=S_{n+1}^{-1}-\frac{\frac{1}{h^2}S_{n+1}^{-1}\mathbf{e}_1 \mathbf{v}_h^{\mathrm{T}}S_{n+1}^{-1}}{1+\frac{1}{h^2}\mathbf{v}_h^{\mathrm{T}}S_{n+1}^{-1}\mathbf{e}_1}\nonumber\\
&=S_{n+1}^{-1}-R_{n+1}.
\end{align}

Our goal is to estimate quite accurately $\|A_h^{-1}\|_p$ with $p\in [1,\infty]$ and with $\|\cdot\|_p$ being the matrix norm induced by the vector norm $\| \mathbf{x}\|_p=\left[\sum |x_j|^p\right]^{1/p}$. We concentrate our efforts in the case where $p=1,2,\infty$, since the other estimates can be obtained via classical interpolation techniques.

We start by estimating $\|S_{n+1}^{-1}\|_1, \|S_{n+1}^{-1}\|_\infty,\|R_{n+1}\|_1$, $\|R_{n+1}\|_\infty$. The latter are used for giving quite precise bounds on  $\|A_h^{-1}\|_1$ and $\|A_h^{-1}\|_\infty$. The estimate for $\|A_h^{-1}\|_2$ can be obtained by a direct check, but it essentially follows from the estimates on $\|A_h^{-1}\|_1$ and $\|A_h^{-1}\|_\infty$, by means of the inequality
 $\|A_h^{-1}\|_2\leq\sqrt{\|A_h^{-1}\|_1\|A_h^{-1}\|_\infty}$.

\subsection{Estimating $\|S_{n+1}^{-1}\|_p$ with $p=1,\infty$}
We have $S_{n+1}^{-1}=h^2T_{n+1}^{-1}$ where $T_{n+1}=T_{n+1}(2-2\cos(\theta))$ and the inverse $\left(T_{n+1}^{-1}\right)_{r,c}=t_r^{(c)}$
\begin{align}
T_{n+1}^{-1}=
\begin{bmatrix*}
t_1^{(1)}&t_1^{(2)}&\ldots&t_1^{(n+1)}\\
t_2^{(1)}&t_2^{(2)}&\ldots&t_2^{(n+1)}\\
\vdots&\vdots&\ddots&\vdots\\
t_{n+1}^{(1)}&t_{n+1}^{(2)}&\ldots&t_{n+1}^{(n+1)}
\end{bmatrix*}=
\begin{bmatrix*}
\mathbf{t}^{(1)}&\mathbf{t}^{(2)}&\ldots&\mathbf{t}^{(n+1)}
\end{bmatrix*}.\nonumber
\end{align}
The components of the inverse $T_{n+1}^{-1}$, $t_r^{(c)}$, are defined by, for a fixed column $c$,
\begin{align}
t_r^{(c)}&=\frac{(n+2-c)r}{n+2},\quad r=1,\ldots, c-1, \text{ for }c>1,\label{eq:tc1}\\
t_r^{(c)}&=\frac{(n+2-r)c}{n+2}, \quad r=c\ldots, n+1,\label{eq:tc2}
\end{align}
and symmetrically for a fixed row $r$
\begin{align}
t_r^{(c)}&=\frac{(n+2-r)c}{n+2},\quad c=1,\ldots, r-1, \text{ for }r>1,\label{eq:tr1}\\
t_r^{(c)}&=\frac{(n+2-c)r}{n+2}, \quad c=r,\ldots, n+1.\label{eq:tr2}
\end{align}

All terms of $S_{n+1}^{-1}$ (and $T_{n+1}^{-1}$) are positive and real, and they are symmetric. Hence by using the explicit expressions of the considered norms, we find
\begin{align}
\|S_{n+1}^{-1}\|_\infty &=\max_r\left\{\sum_{c=1}^{n+1}\left(S_{n+1}^{-1}\right)_{r,c}\right\}=\max_r\left\{h^2\sum_{c=1}^{n+1}\left(T_{n+1}^{-1}\right)_{r,c}\right\}\nonumber \\
&=\max_c\left\{h^2\sum_{r=1}^{n+1}\left(T_{n+1}^{-1}\right)_{r,c}\right\}=\max_c\left\{\sum_{r=1}^{n+1}\left(S_{n+1}^{-1}\right)_{r,c}\right\}=\|S_{n+1}^{-1}\|_1.
\end{align}

Numerically it is obvious that the highest row sum for matrices $T_{n+1}^{-1}$ with $n+1$ even is for row index $r=(n+1)/2$ (or $r=(n+1)/2+1$, they are equal). For odd $n+1$, the highest row sum is for row index $r=(n+2)/2$.

Thus for $n+1$ even
\begin{align}
\|T_{n+1}^{-1}\|_\infty&=\sum_{c=1}^{n+1} t_{(n+1)/2}^{(c)}=\sum_{c=1}^{(n-1)/2}\frac{(n+2-(n+1)/2)c}{n+2}+\sum_{c=(n+1)/2}^{n+1}\frac{(n+2-c)(n+1)/2}{n+2}\nonumber\\
&=\frac{(n+1)^2+2(n+1)}{8}=\frac{1+2h}{8h^2}\nonumber
\end{align}
and for $n+1$ odd
\begin{align}
\|T_{n+1}^{-1}\|_\infty&=\sum_{c=1}^{n+1} t_{(n+2)/2}^{(c)}=\sum_{c=1}^{n/2}\frac{(n+2-(n+2)/2)c}{n+2}+\hspace{-1em}\sum_{c=(n+2)/2}^{n+1}\hspace{-1em}\frac{(n+2-c)(n+2)/2}{n+2}\nonumber\\
&=\frac{(n+2)^2}{8}=\frac{(1+h)^2}{8h^2}\nonumber
\end{align}
Consequently, for $n+1$ even, we deduce
\begin{align}
\|S_{n+1}^{-1}\|_\infty&=h^2\|T_{n+1}^{-1}\|_\infty=\frac{1+2h}{8},\label{eq:Sne}
\end{align}
and for $n+1$ odd, we have
\begin{align}
\|S_{n+1}^{-1}\|_\infty&=h^2\|T_{n+1}^{-1}\|_\infty=\frac{1+2h+h^2}{8}.\label{eq:Sno}
\end{align}
As a conclusion, for all $n+1$ and using the symmetry and \eqref{eq:Sne} and \eqref{eq:Sno}, we obtain that
\begin{align}
\|S_{n+1}^{-1}\|_\infty=\|S_{n+1}^{-1}\|_1\leq\frac{1+2h+h^2}{8},
\end{align}
and the  limit as the matrix size tends to infinity, that is $h\to0$, is $\|S_{n+1}^{-1}\|_\infty=\|S_{n+1}^{-1}\|_1\to\frac{1}{8}$.

\subsection{Estimating $\|R_{n+1}\|_p$  for $p=1,\infty$}
Since $S_{n+1}^{-1}=h^2T_{n+1}^{-1}$ and $T_{n+1}^{-1}\mathbf{e}_1=\mathbf{t}^{(1)}$, we find that
\begin{align}
R_{n+1}&=\frac{\frac{1}{h^2}S_{n+1}^{-1}\mathbf{e}_1\mathbf{v}_h^{\mathrm{T}}S_{n+1}^{-1}}{1+\frac{1}{h^2}\mathbf{v}_h^{\mathrm{T}}S_{n+1}^{-1}\mathbf{e}_1}\nonumber\\
&=
\frac{T_{n+1}^{-1}\mathbf{e}_1\mathbf{v}_h^{\mathrm{T}}S_{n+1}^{-1}}{1+\mathbf{v}_h^{\mathrm{T}}T_{n+1}^{-1}\mathbf{e}_1}\nonumber\\
&=\frac{\mathbf{t}^{(1)}\mathbf{v}_h^{\mathrm{T}}S_{n+1}^{-1}}{1+\mathbf{v}_h^{\mathrm{T}}\mathbf{t}^{(1)}}.\label{eqRn}
\end{align}
Moreover we have from \eqref{eq:tc2} that the components of $\mathbf{t}^{(1)}$ are
\begin{align}
t_r^{(1)}&=\frac{n+2-r}{n+2}=1-\frac{hr}{1+h}=\frac{1-h(r-1)}{1+h}, \quad r=1,\ldots, n+1.\label{eqtr1}
\end{align}
and we have from \eqref{eq:Tn}
\begin{align}
\mathbf{v}_h^{\mathrm{T}}&=\begin{bmatrix*}
\vartheta h^2-2&(1-\vartheta)h^2+1&0&\ldots &0
\end{bmatrix*}=\begin{bmatrix*}v_1&v_2&0&\ldots&0\end{bmatrix*}.\nonumber
\end{align}
Thus
\begin{align}
\mathbf{v}_h^{\mathrm{T}}\mathbf{t}^{(1)}&=v_1t_1^{(1)}+v_2t_2^{(1)}
=(\vartheta h^2-2)\frac{1}{1+h}+((1-\vartheta)h^2+1)\frac{1-h}{1+h}\nonumber\\
&=\frac{(\vartheta-1) h^3 +h^2-h-1}{1+h},\nonumber
\end{align}
and
\begin{align}
1+\mathbf{v}_h^{\mathrm{T}}\mathbf{t}^{(1)}&=\frac{h^2\left((\vartheta-1)h+1\right)}{1+h}.\label{eq:1pvt1}
\end{align}
Also
\begin{align}
\mathbf{v}_h^{\mathrm{T}}S_{n+1}^{-1} & =h^2\mathbf{v}_h^{\mathrm{T}}T_{n+1}^{-1}\nonumber\\
&=h^2\begin{bmatrix*}
v_1t_1^{(1)}+v_2t_2^{(1)}&v_1t_1^{(2)}+v_2t_2^{(2)}&\ldots&v_1t_1^{(c)}+v_2t_2^{(c)}&\ldots&v_1t_1^{(n+1)}+v_2t_2^{(n+1)}
\end{bmatrix*}.\nonumber
\end{align}
and thus the components of the row vector $\left(\mathbf{v}_h^{\mathrm{T}}S_{n+1}^{-1}\right)_c$ are
\begin{align}
\left(\mathbf{v}_h^{\mathrm{T}}S_{n+1}^{-1}\right)_c&=h^2\left(v_1t_1^{(c)}+v_2t_2^{(c)}\right),\nonumber
\end{align}
and the components of the matrix $\left(\mathbf{t}^{(1)}\mathbf{v}_h^{\mathrm{T}}S_{n+1}^{-1}\right)_{r,c}$ are
\begin{align}
\left(\mathbf{t}^{(1)}\mathbf{v}_h^{\mathrm{T}}S_{n+1}^{-1}\right)_{r,c}&=t_r^{(1)}h^2\left(v_1t_1^{(c)}+v_2t_2^{(c)}\right),
\end{align}
where $t_r^{(1)}$ is defined in \eqref{eqtr1}, and
\begin{align}
t_1^{(c)}&=\frac{1-h(c-1)}{1+h}, \quad c=1,\ldots, n+1,\\
t_2^{(c)}&=\begin{cases*}
\frac{1-h}{1+h},&$c=1$,\\
\frac{2-2h(c-1)}{1+h},&$c >1$,
\end{cases*}
\end{align}
are defined in \eqref{eq:tr1} and \eqref{eq:tr2}.
Therefore for $c=1$
\begin{align}
\left(\mathbf{t}^{(1)}\mathbf{v}_h^{\mathrm{T}}S_{n+1}^{-1}\right)_{r,1}&=t_r^{(1)}h^2\left(v_1t_1^{(1)}+v_2t_2^{(1)}\right)=t_r^{(1)}h^2\mathbf{v}_h^{\mathrm{T}}\mathbf{t}^{(1)}\nonumber\\
&=\frac{1-h(r-1)}{1+h}h^2\frac{(\vartheta-1)h^3+h^2-h-1}{1+h},\nonumber\\
&=\frac{h^2 (1 - h (r - 1)) ((\vartheta - 1)h^3+ h^2 - h - 1)}{(h + 1)^2},\label{eq:tvsn1}
\end{align}
and for  $c>1$
\begin{align}
\left(\mathbf{t}^{(1)}\mathbf{v}_h^{\mathrm{T}}S_{n+1}^{-1}\right)_{r,c}&=t_r^{(1)}h^2\left(v_1t_1^{(c)}+v_2t_2^{(c)}\right)=t_r^{(1)}h^2\mathbf{v}_h^{\mathrm{T}}\mathbf{t}^{(c)}\nonumber\\
&=\frac{1-h(r-1)}{1+h}h^2\left(\left(\vartheta h^2-2\right)\frac{1-h(c-1)}{1+h}  +\left((1-\vartheta)h^2+1\right)\frac{2-2h(c-1)}{1+h}\right)\nonumber\\
&=\frac{h^4 (2-\vartheta ) (1-h(c -1) ) (1-h(r-1))}{(h + 1)^2}\label{eq:tvsn2}
\end{align}
Thus we can now define the components of $(R_{n+1})_{r,c}=\left(\frac{\mathbf{t}^{(1)}\mathbf{v}_h^{\mathrm{T}}S_{n+1}^{-1}}{1+\mathbf{v}_h^{\mathrm{T}}\mathbf{t}^{(1)}}\right)_{r,c}$, defined in \eqref{eqRn}, since we have \eqref{eq:1pvt1}, \eqref{eq:tvsn1}, and \eqref{eq:tvsn2}. For $c=1$ we have
\begin{align}
\left(R_{n+1}\right)_{r,1}&=\frac{h^2 (1 - h (r - 1)) ((\vartheta - 1)h^3+ h^2 - h - 1)}{(h + 1)^2} \bigg/ \frac{h^2\left((\vartheta-1)h+1\right)}{1+h}\nonumber\\
&=\frac{h(r-1)-1}{h(\vartheta-1)+1}-\frac{h^2(h(r-1)-1)}{h+1}\label{eq:rnr1}
\end{align}
and for $c>1$
\begin{align}
\left(R_{n+1}\right)_{r,c}&=\nonumber\\
&=\frac{h^4 (2-\vartheta ) (1-h(c -1) ) (1-h(r-1))}{(h + 1)^2}\bigg/ \frac{h^2\left((\vartheta-1)h+1\right)}{1+h}\nonumber\\
&=\frac{h^2 (2 - \vartheta) (1 - (c - 1) h) (1 - h (r - 1))}{(h + 1) (h (\vartheta- 1) + 1)}\label{eq:rnrc}
\end{align}
Numerically it is obvious that $\|R_{n+1}\|_1$ and $\|R_{n+1}\|_\infty$ is always for the first column and first row.

Now we compute $\|R_{n+1}\|_1$,
\begin{align}
\|R_{n+1}\|_1&=\sum_{r=1}^{n+1}\left|\left(R_{n+1}\right)_{r,1}\right|\nonumber\\
&=\sum_{r=1}^{n+1}\left|\frac{h(r-1)-1}{h(\vartheta-1)+1}-\frac{h^2(h(r-1)-1)}{h+1}\right|\nonumber\\
&=\sum_{r=1}^{n+1}\left|\frac{(1 - h (r - 1)) (h^3 (\vartheta - 1) + h^2 - h - 1)}{(h + 1) (h (\vartheta - 1) + 1)}\right|\nonumber\\
&=\frac{-h^3 (\vartheta - 1) - h^2 + h +1}{(h + 1) (h (\vartheta - 1) + 1)}\sum_{r=1}^{n+1}(1 +h - hr ))\nonumber\\
&=\frac{-h^3 (\vartheta - 1) - h^2 + h +1}{(h + 1) (h (\vartheta - 1) + 1)}\left((n+1)(1+h) - h \frac{(n+1)(n+2)}{2}\right)\nonumber\\
&=\frac{-h^3 (\vartheta - 1) - h^2 + h +1}{(h + 1) (h (\vartheta - 1) + 1)}\left(\frac{1+h}{h} - \frac{1+h}{2h}\right)\nonumber\\
&=\frac{h^3 (1-\vartheta ) - h^2 + h +1}{2h (h (\vartheta - 1) + 1)},\label{eq:rn1}
\end{align}
since $0<\vartheta<1$.

We now compute $\|R_{n+1}\|_\infty$, by taking into account that all coefficients are positive except the first in the first column. We have from \eqref{eq:rnr1} and \eqref{eq:rnrc}
\begin{align}
\left(R_{n+1}\right)_{1,1}&=\frac{h^2}{h+1}-\frac{1}{h(\vartheta-1)+1}\label{eqrna}\\
\left(R_{n+1}\right)_{1,c}&=\frac{h^2 (2 - \vartheta) (1 - (c - 1) h)}{(h + 1) (h (\vartheta- 1) + 1)},\quad c=2,\ldots,n+1\label{eqrnb}
\end{align}
Thus
\begin{align}
\|R_{n+1}\|_\infty&=-\left(R_{n+1}\right)_{1,1}+\sum_{c=2}^{n+1}\left(R_{n+1}\right)_{1,c}\nonumber\\
&=-\left(\frac{h^2}{h+1}-\frac{1}{h(\vartheta-1)+1}\right)+\sum_{c=2}^{n+1}\frac{h^2 (2 - \vartheta) (1 - (c - 1) h)}{(h + 1) (h (\vartheta- 1) + 1)}\nonumber\\
&=-\left(\frac{h^2}{h+1}-\frac{1}{h(\vartheta-1)+1}\right)+\frac{h^2 (2 - \vartheta) }{(h + 1) (h (\vartheta- 1) + 1)}\sum_{c=2}^{n+1}(1+h - c h)\nonumber\\
&=-\left(\frac{h^2}{h+1}-\frac{1}{h(\vartheta-1)+1}\right)+\frac{h^2 (2 - \vartheta) }{(h + 1) (h (\vartheta- 1) + 1)}\left((n+1-1)(1+h) - h\left(\frac{(n+1)(n+2)}{2}-1\right)\right)\nonumber\\
&=-\left(\frac{h^2}{h+1}-\frac{1}{h(\vartheta-1)+1}\right)+\frac{h^2 (2 - \vartheta) }{(h + 1) (h (\vartheta- 1) + 1)}\left(\frac{1-h^2}{h} - \frac{1+h-2h^2}{2h}\right)\nonumber\\
&=-\left(\frac{h^2}{h+1}-\frac{1}{h(\vartheta-1)+1}\right)+\frac{1}{2}\frac{h (2 - \vartheta) (1-h)}{(h + 1) (h (\vartheta- 1) + 1)}\nonumber\\
&=\frac{1}{2}\frac{h (2 - \vartheta) (1-h)-2h^2(h (\vartheta- 1) + 1)+2(h+1)}{(h + 1) (h (\vartheta- 1) + 1)}\nonumber\\
&\to 1 \text{ as } h\to 0.\nonumber
\end{align}

\subsection{Estimating $\|A_h^{-1}\|_p$  for $p=1,\infty$}
Numerically it is obvious that $\|A_h^{-1}\|_1$ is computed on the first column, thus since $A_h^{-1}=S_{n+1}^{-1}-R_{n+1}$,  $S_{n+1}^{-1}$ and $R_{n+1}$ positive and negative respectively, we can just compute the norm directly for $A_h^{-1}$. The sum of the positive elements of the first column of $T_{n+1}^{-1}$ is equal to $\frac{{n+1}}{2}$, and thus the sum for the first column of $S_{n+1}^{-1}$ is $\frac{h^2(n+1)}{2}=\frac{h}{2}$, and the sum of components $-(R_{n+1})_{r,1}$ is given in \eqref{eq:rn1}, that is
\begin{align}
\|A_h^{-1}\|_1&=\frac{h}{2}+\frac{h^3 (1-\vartheta ) - h^2 + h +1}{2h (h (\vartheta - 1) + 1)}\nonumber\\
&=\frac{ h^3 (\vartheta - 1) + h^2+h^3 (1-\vartheta ) - h^2 + h +1}{2h (h (\vartheta - 1) + 1)}\nonumber\\
&=\frac{h+1}{2h (h (\vartheta - 1) + 1)}\label{eq:normah1}
\end{align}

Now we compute $\|A_h^{-1}\|_\infty$. $A_h^{-1}=S_{n+1}^{-1}-R_{n+1}$. We have from \eqref{eq:tr2}
\begin{align}
\left(S_{n+1}^{-1}\right)_{1,c}&=h^2t_1^{(c)}=h^2\frac{1+h(1-c)}{1+h}, \quad c=1,\ldots, n+1\nonumber
\end{align}
and by using \eqref{eqrna}  we get
\begin{align}
\left(A_h^{-1}\right)_{1,1}&=(S_{n+1}^{-1})_{1,1}-\left(R_{n+1}\right)_{1,1}\nonumber\\
&=\frac{h^2}{1+h}-\left(\frac{h^2}{h+1}-\frac{1}{h(\vartheta-1)+1}\right)\nonumber\\
&=\frac{1}{h(\vartheta-1)+1}\label{eq:ahi11}
\end{align}
and by \eqref{eqrnb}
\begin{align}
\left(A_h^{-1}\right)_{1,c}&=(S_{n+1}^{-1})_{1,c}-\left(R_{n+1}\right)_{1,c}\nonumber\\
&=h^2\frac{1+h(1-c)}{1+h}-\frac{h^2 (2 - \vartheta) (1 - (c - 1) h)}{(h + 1) (h (\vartheta- 1) + 1)}\nonumber\\
&=\frac{h^2(1+h(1-c))}{1+h}\left(1-\frac{ 2 - \vartheta}{h (\vartheta- 1) + 1}\right)\nonumber\\
&=\frac{h^2(1+h(1-c))( \vartheta-1)}{h (\vartheta- 1) + 1}.\label{eq:ahi1c}
\end{align}
Since $\left(A_h^{-1}\right)_{1,1}$ in \eqref{eq:ahi11} is always positive and $\left(A_h^{-1}\right)_{1,c}$ of \eqref{eq:ahi1c} is always negative we have
\begin{align}
\|A_h^{-1}\|_\infty&=\left(A_h^{-1}\right)_{1,1}-\sum_{c=2}^{n+1}\left(A_h^{-1}\right)_{1,c}\nonumber\\
&=\frac{1}{h(\vartheta-1)+1}-\sum_{c=2}^{n+1}\frac{h^2(1+h(1-c))( \vartheta-1)}{h (\vartheta- 1) + 1}\nonumber\\
&=\frac{1}{h(\vartheta-1)+1}-\frac{h^2( \vartheta-1)}{h (\vartheta- 1) + 1}\sum_{c=2}^{n+1}(1+h-hc)\nonumber\\
&=\frac{1}{h(\vartheta-1)+1}-\frac{h^2( \vartheta-1)}{h (\vartheta- 1) + 1}\left((n+1-1)(1+h)-h\left(\frac{(n+1)(n+2)}{2}-1\right)\right)\nonumber\\
&=\frac{1}{h(\vartheta-1)+1}-\frac{h^2( \vartheta-1)}{h (\vartheta- 1) + 1}\left(\frac{1-h^2}{h}-\frac{1+h-2h^2}{2h}\right)\nonumber\\
&=\frac{2-h( \vartheta-1)(1-h)}{2(h(\vartheta-1)+1)}\label{eq:normahinf}
\end{align}

As a conclusion we deduce from \eqref{eq:normah1} and \eqref{eq:normahinf}
\begin{align}
\|A_h^{-1}\|_2&\leq\sqrt{\|A_h^{-1}\|_1\|A_h^{-1}\|_\infty}\nonumber\\
&=\sqrt{\frac{h+1}{2h (h (\vartheta - 1) + 1)}\frac{2-h( \vartheta-1)(1-h)}{2(h(\vartheta-1)+1)}}\nonumber\\
&=\frac{1}{2(h(\vartheta-1)+1)}\sqrt{\frac{h+1}{h}(2-h( \vartheta-1)(1-h))}\nonumber\\
&=\frac{1}{2(h(\vartheta-1)+1)}\sqrt{\frac{2(h+1)+(h^2-1)h(\vartheta-1)}{h}}\nonumber\\
&=\frac{1}{2(h(\vartheta-1)+1)}\sqrt{\frac{2}{h}+2+(h^2-1)(\vartheta-1)}\nonumber
\end{align}

In order to make a comparison, we recall that we know the exact asymptotical behavior of $\|S_{n+1}^{-1}\|_2$, with $S_{n+1}$ being the pure Toeplitz counterpart of $A_h$, as reported below
\begin{align}
\|S_{n+1}^{-1}\|_2&=\frac{1}{\lambda_{\min}(S_{n+1})}=\frac{h^2}{4\sin^2\left(\frac{\pi}{2(n+2)}\right)}=\left(\frac{h}{2\sin\left(\frac{\pi h}{2(1+h)}\right)}\right)^2\stackrel{h\to 0}{\to } \frac{1}{\pi^2}.
\end{align}

\subsection{Spectral results: comments}\label{ssec:spectral 1D}

Here we give a short discussion on few items that, for some aspects, will be considered in more detail in Section \ref{sec:2D} and for other aspects will be listed as open problems in the conclusion Section  \ref{sec:final}.

\begin{itemize}
\item The estimates for $\|A_h^{-1}\|_p$ are tight and the growth is like $n^{1/p}$: however the numerical growth of the error seems to be bounded by a constant independently of $p$. The reason relies on the vectors for which the norm is attained. Such vectors should be concentrated on the first component and this is quite unphysical and it is not observed in practice.
\item Even if $A_h$ and its inverse are not symmetric we can prove the spectrum of the related matrix-sequence is clustered along a real positive interval, using the results of the GLT technology reported in Subsection \ref{ssec:glt} (see also  \cite{BaSe,GoSe}): we refer to Subsection \ref{ssec:spectral 2D} where the analysis is performed both in $1\mathbb{D}$ and $2\mathbb{D}$.
\item Regarding the estimates of $\|A_h^{-1}\|_p$,  the $2\mathbb{D}$ case (and generically the $d\mathbb{D}$ case) is more difficult, but we can take advantage of the one dimensional case and from a clever tensor structure of the problem when the domain is rectangular (hyper-rectangular in the $d\mathbb{D}$ case).
\item When the domain is generic a possibility is given by embedding techniques already exploited in the distributional setting via the GLT approach (see \cite{GLT1,GLT2}).
\end{itemize}


\subsection{Numerical tests in $1\mathbb{D}$}\label{ssec:numerical_tests}
We consider the $1\mathbb{D}$ problem~\eqref{mainProblem}
with $a=0$ and $b=\pi$. We choose $f=-\sin(x)$, $g_a=0$, and $g_b=0$ so that $\mathbf{u}=-\sin(x_i)$ is the exact solution in points $x_i$.

We perform several tests varying the value of $\vartheta \in [0,1]$, in order to establish whether the convergence of the method depends on the choice of $\vartheta$.
In practice, we choose  $\vartheta$ and $n$ and we compute $h$ and $x_0$ accordingly:
\[
h = \frac{b-a}{n+\vartheta}, \qquad x_0 = b - (n+1)h.
\]

\begin{figure}
\centering
\includegraphics[width=0.99\textwidth]{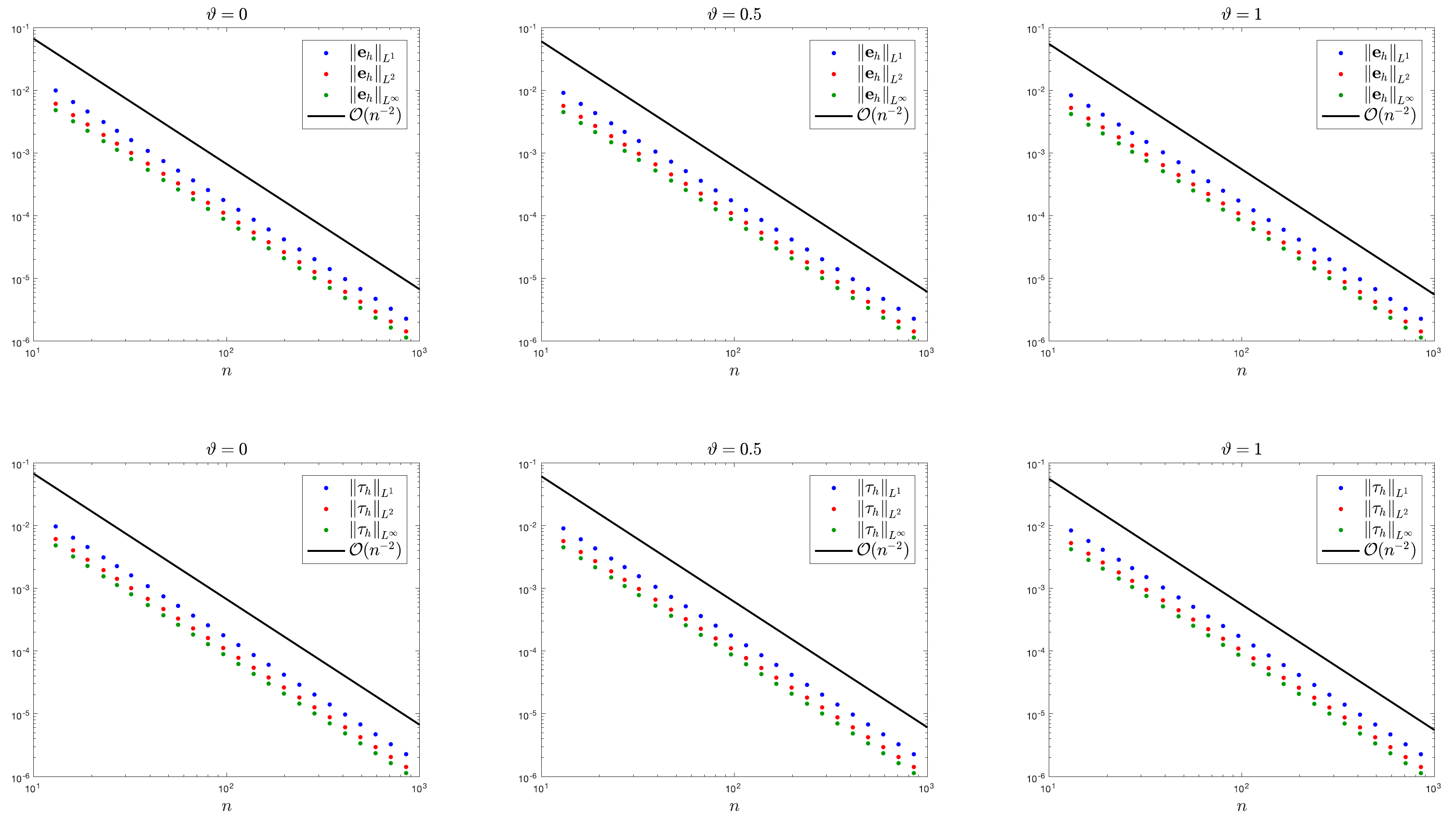}
\caption{The dots represent the $p$ - norm of the numerical error $\mathbf{e}_h$ (top) and consistency error $\tau_h$ (bottom) for different values of $n$ (horizontal axis) and $\vartheta$:
$\vartheta=0$ (left), $\vartheta=0.5$ (middle), $\vartheta=1$ (right). The solid line is a reference for second-order decay. }
\label{fig:second_order}
\end{figure}

The numerical error $\mathbf{e}_h = \mathbf{u} - \mathbf{u}_h$ satisfies the following equation:
\begin{equation*}
A_h \mathbf{e}_h = \tau_h,
\end{equation*}
where $\tau_h$ is the consistency error:
\[
\tau_h = \mathbf{f}_h - A_h \mathbf{u}.
\]
Consider the $p-$norm:
\begin{equation} \left\|\tau_h\right\|_{L^p} \approx \left(h\sum_{i=0}^{n} \left|\tau_h(x_i)\right|^p \right)^{\frac{1}{p}},
\end{equation}
\begin{equation} \left\| \mathbf{e}_h \right\|_{L^p} \approx \left(h  \sum_{i=0}^{n} \left| \mathbf{e}_h(x_i)\right|^p \right)^{\frac{1}{p}},
\end{equation}
In Fig.~\ref{fig:second_order} we show that:
\begin{equation} \left\|\tau_h \right\|_{L^p}, \left\|\mathbf{e}_h \right\|_{L^p} \approx O(h^2), \ \textrm{for} \ p=1, 2, \infty,
\end{equation}
confirming that the method is second-order consistent and accurate.

\begin{figure}[htp]
\centering
\includegraphics[width=0.48\textwidth]{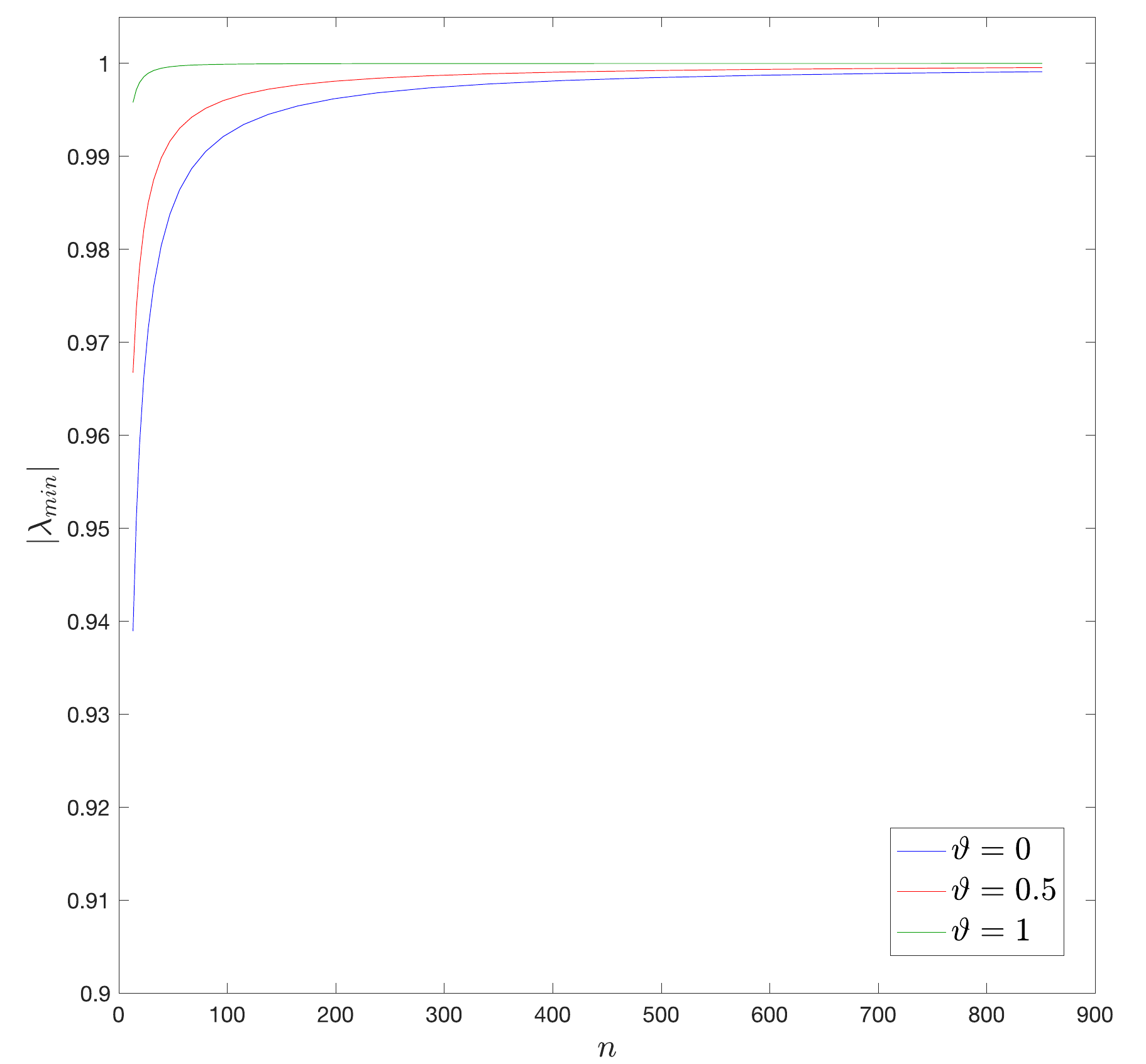}
\includegraphics[width=0.48\textwidth]{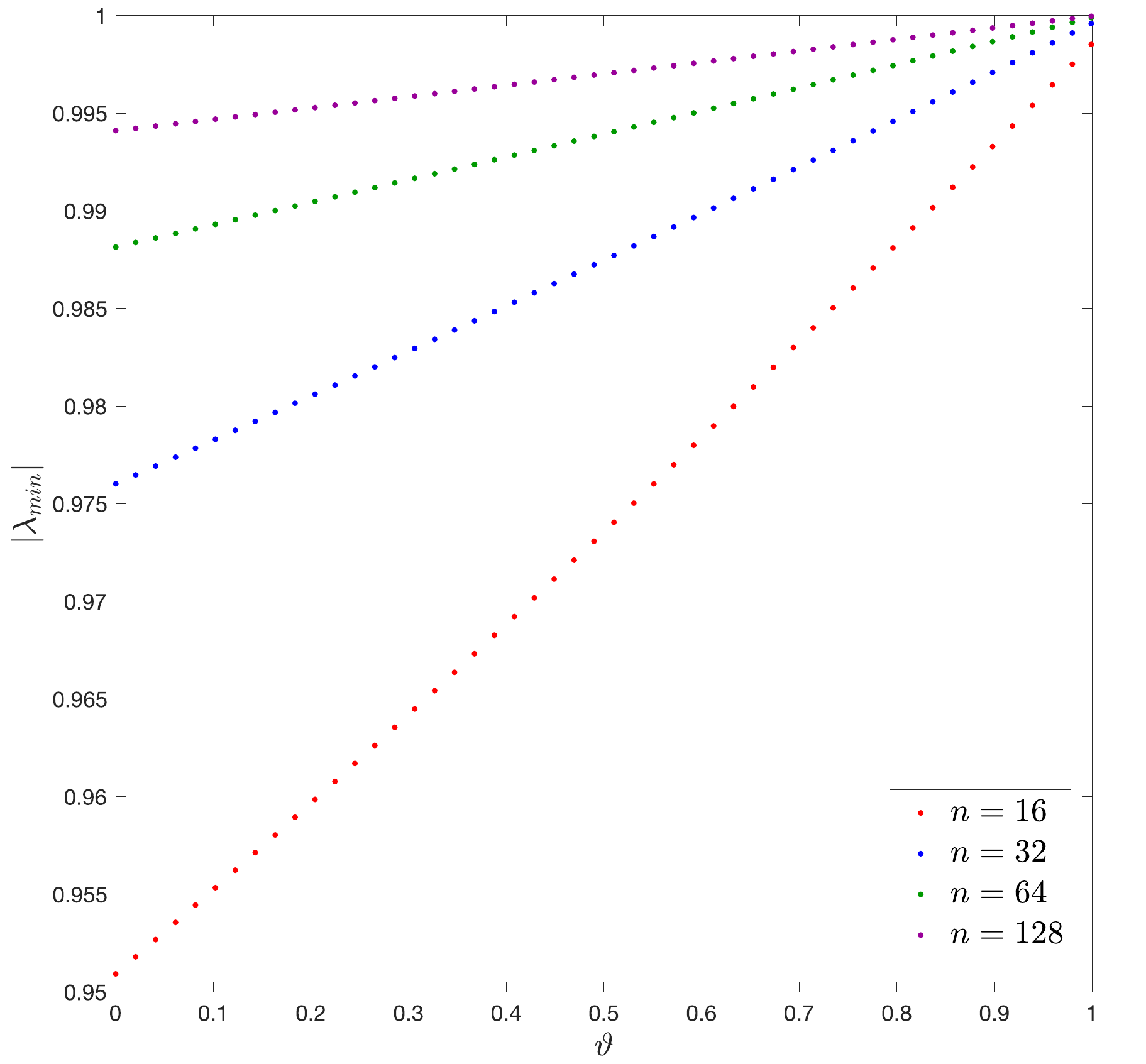}
\label{fig:lambdamin}
\caption{Smallest eigenvalue in absolute value (vertical axis) for different values of $n$ (horizontal axis, left plot) or for different values of $\vartheta$ (horizontal axis, right plot).}
\end{figure}

We complete the analysis showing the behaviour of the spectral radius of the matrix $A_h^{-1}$. Fig.~\ref{fig:lambdamin} shows how the smallest eigenvalue (in absolute value) of the matrix $A_h$ changes in relation to $n$ (left panel) and in relation to $\vartheta$ (right panel).

Fig.~\ref{fig:lambdamin} shows that the smallest eigenvalue (in absolute value) of the matrix $A_h$ essentially does not depend on the value of $\vartheta$ and it approaches a constant values when $n$ goes to infinity.\\

Since $\left\|\mathbf{e}_h\right\|_{L^p} \leq \left\|A_h^{-1}\right\|_{p} \left\|\tau_h\right\|_{L^p}$, we can conclude that $\left\|\mathbf{e}_h\right\|_{L^p} \approx O(h^2)$ and $\left\|A_h^{-1}\right\|_{p} \left\|\tau_h\right\|_{L^p} \approx O(h^{2 - \frac{1}{p}})$, as predicted in the first item of Subsection~\ref{ssec:spectral 1D}.

\section{Problem formulation in $2\mathbb{D}$ and related analysis}\label{sec:2D}

The section is organized into three parts: first we introduce the $d$-level notation and the $d$-level Toeplitz matrices in Subsection \ref{ssec:multilevel}, secondly we define the notion of spectral and singular value distribution and the $*$-algebra of Generalized Locally Toeplitz matrix-sequences in Subsection \ref{ssec:glt}, then we describe the matrices arising in the approximation of a Dirichlet problem by the Coco--Russo method in Subsection \ref{ssec:method 2D}, and finally we give a spectral analysis of the resulting matrix-sequences in Subsection \ref{ssec:spectral 2D}.

\subsection{Multilevel notation: the case of multilevel Toeplitz and diagonal sampling matrices}
	\label{ssec:multilevel}

We start by introducing the multi-index notation, which is useful in our context.
A multi-index $\mathbf i\in\mathbb Z^d$, also called a $d$-index, is simply a (row) vector in $\mathbb Z^d$; its components are denoted by $i_1,\ldots,i_d$.

\begin{itemize}
	\item $\mathbf0,\,\mathbf1,\,\mathbf2,\,\ldots$ are the vectors of all zeros, all ones, all twos, $\ldots$ (their size will be clear from the context).
	\item For any $d$-index $\mathbf m$, we set $N(\mathbf m)=\prod_{j=1}^dm_j$ and we write $\mathbf m\to\infty$ to indicate that $\min(\mathbf m)\to\infty$.
  \item If $\mathbf h, \mathbf k$ are $d$-indices, $\mathbf h\le\mathbf k$ means that $h_r\le k_r$ for all $r=1,\ldots,d$.
  \item The standard lexicographic ordering is assumed uniformly
\begin{equation}
\left[\ \ldots\ \left[\ \left[\ (j_1,\ldots,j_d)\ \right]_{j_d=h_d,\ldots,k_d}\ \right]_{j_{d-1}=h_{d-1},\ldots,k_{d-1}}\ \ldots\ \right]_{j_1=h_1,\ldots,k_1}.
\label{eq:standardLex}
\end{equation}
\end{itemize}
For instance, in the case $d=2$ the ordering is the following: $(h_1,h_2)$, $(h_1,h_2+1)$, $\ldots,$ $(h_1,k_2)$, $(h_1+1,h_2)$,

\noindent{\bf Multilevel Toeplitz Matrices.}

We now briefly summarize the definition and few relevant properties of multilevel Toeplitz matrices, that we will employ in the analysis of the $2\mathbb{D}$ setting.
Given $\mathbf{n}\in \mathbb{N}^d$, a matrix of the form
\begin{equation*}
  [a_{\mathbf{i}-\mathbf{j}}]_{\mathbf{i},\mathbf{j}=\mathbf{e}}^{\mathbf{n}} \in \mathbb{C}^{N(\mathbf{n}) \times N(\mathbf{n})}
\end{equation*}
with $\mathbf{e}$ vector of all ones, with entries $a_\mathbf{k}
\in \mathbb{C}$, $\mathbf{k} =
-(\mathbf{n}-\mathbf{e}), \ldots, \mathbf{n}-\mathbf{e}$, is
called a multilevel Toeplitz matrix, or, more precisely, a
$d$-level Toeplitz matrix. Let $\phi : [-\pi, \pi]^{d}
\rightarrow \mathbb{C}^{r\times r}$ a matrix-valued function in which each entry belongs to $L^1([-\pi, \pi]^d)$.
 We denote the Fourier coefficients of the generating function $\phi$ as
\begin{equation*}
  \hat{\phi}_\mathbf{k} =\frac{1}{(2\pi)^d} \int_{[-\pi,\pi]^d}
  \phi(\boldsymbol{\theta})e^{-\hat i (\mathbf{k}, \boldsymbol{\theta})}
  \ d{\boldsymbol \theta} \in \mathbb{C}, \quad \mathbf{k} \in \mathbb{Z}^d,
\end{equation*}
where the integrals are computed component-wise and $(\mathbf{k},
\boldsymbol{\theta}) = k_1\theta_1 + \ldots + k_d\theta_d$. For
every $\mathbf{n} \in \mathbb{N}^d$, the $\mathbf{n}$-th Toeplitz
matrix associated with $\phi$ is defined as
\begin{equation*}
T_\mathbf{n}(\phi) :=
[\hat{\phi}_{\mathbf{i}-\mathbf{j}}]_{\mathbf{i},\mathbf{j}=\mathbf{e}}^{\mathbf{n}}
\end{equation*}
or, equivalently, as
\begin{equation}\label{d-level T}
T_\mathbf{n}(\phi) = \sum_{|j_1|<n_1} \ldots \sum_{|j_d|<n_d}
\hat{\phi}_{(j_1,\ldots, j_d)}[J_{n_1}^{(j_1)} \otimes \ldots \otimes J_{n_d}^{(j_d)}],
\end{equation}
where $\otimes$ denotes the (Kronecker) tensor product of
matrices, while $J_m^{(l)}$ is the matrix of order $m$ whose
$(i,j)$ entry equals $1$ if $i-j=l$ and zero otherwise.
We call $\{T_\mathbf{n}(\phi)\}_{\mathbf{n}\in \mathbb{N}^d}$ the
family of (multilevel block) Toeplitz matrices associated with
$\phi$, which, in turn, is called the generating function of
$\{T_\mathbf{n}(\phi)\}_{\mathbf{n}\in \mathbb{N}^d}$ .

\noindent{\bf Multilevel Diagonal Sampling Matrices.}
For $n\in\mathbb N$ and $a:[0,1]\to\mathbb C$, we define the diagonal sampling matrix $D_n(a)$ as the diagonal matrix
\begin{equation*}
  D_n(a)=\mathop{\rm diag}_{i=1,\ldots,n}a\Bigl(\frac{i}{n}\Bigr)=\begin{bmatrix} a(\frac1n) & & & \\ & a(\frac2n) & & \\ & & \ddots & \\
& & & a(1)\end{bmatrix} \in\mathbb C^{n\times n}.
\end{equation*}
For $\mathbf{n}\in\mathbb N^d$ and $a:[0,1]^d\to\mathbb C$, we define the multilevel diagonal sampling matrix $D_\mathbf{n}(a)$ as the diagonal matrix
\begin{equation*}
  D_\mathbf{n}(a)=\mathop{\rm diag}_{\mathbf{i}=\mathbf1,\ldots,\mathbf{n}}a\Bigl(\frac{\mathbf{i}}{\mathbf{n}}\Bigr)\in\mathbb C^{N(\mathbf{n})\times N(\mathbf{n})},
\end{equation*}
with the lexicographical ordering~\eqref{eq:standardLex} as discussed at the beginning of the subsection.

\subsection{GLT matrix-sequences: operative features}
	\label{ssec:glt}

We start with the definition of distribution in the sense of the eigenvalues (spectral distribution) and in the sense of the singular values (singular value distribution) for a given matrix-sequence. Then we give the operative feature of the $*$-algebra of matrix-sequences.

\newtheorem{definition}[theorem]{Definition}

\begin{definition}\label{def-distribution}
Let $\{A_n\}_n$ be a sequence of matrices, with $A_n$ of size
$d_n$, and let $f:D\subset\mathbb R^t\to\mathbb{C}$ be
a measurable function defined on a set $D$ with
$0<\mu_t(D)<\infty$.
\begin{itemize}
\item We say that $\{A_n\}_n$ has a (asymptotic) singular value distribution described by $f$, and we write $\{A_n\}_n\sim_\sigma f$, if
  \begin{equation}\label{distribution:sv-sv}
    \lim_{n\to\infty}\frac1{d_n}\sum_{i=1}^{d_n}F(\sigma_i(A_n))=\frac1{\mu_t(D)}\int_D F(|f(\mathbf x)|) \ d{\mathbf x},\quad\forall\,F\in C_c(\mathbb R).
  \end{equation}
\item We say that $\{A_n\}_n$ has a (asymptotic) spectral (or eigenvalue) distribution described by $f$, and we write $\{A_n\}_n\sim_\lambda f$, if
  \begin{equation}\label{distribution:sv-eig}
    \lim_{n\to\infty}\frac1{d_n}\sum_{i=1}^{d_n}F(\lambda_i(A_n))=\frac1{\mu_t(D)}\int_D F(f(\mathbf x))\ d{\mathbf x},\quad\forall\,F\in C_c(\mathbb C).
\end{equation}
\end{itemize}
If $\{A_n\}_n$ has both a singular value and an eigenvalue
distribution described by $f$, then we write $\{A_n\}_n\sim_{\sigma,\lambda}f$.
\end{definition}
The symbol $f$ contains spectral/singular value information briefly described informally as follows. With reference to (\ref{distribution:sv-eig}), assuming that $d_n$ is large enough and $f$ is at least Riemann integrable, except possibly for a small number of outliers, the eigenvalues of $A_n$ are approximately formed by the samples of $f$ over a uniform grid in $D$, so that the range of $f$ is a (weak) cluster for the eigenvalues of $\{A_n\}_n$.
It is then clear that the symbol $f$ provides a `compact' and a quite accurate description of the spectrum of the matrices $A_n$  for $n$ large enough.
Relation (\ref{distribution:sv-sv}) has the same meaning when talking of the singular values of $A_n$ and by replacing
$f$ with $|f|$.
	
	A $d$-level ($d\ge 1$ integer) GLT matrix-sequence $\{A_n\}_n$ is nothing more than a matrix--sequence endowed with a measurable function $\kappa:[0,1]^d\times[-\pi,\pi]^d\to\mathbb C$ called \emph{symbol} characterizing the distributional properties of its singular values, and, under certain hypothesis, of its spectrum. For a complete overview of the theory we refer to the books~\cite{MR3674485,GLTvol2}, while here we recall only the operative features we need for our restricted setting. Since we have already introduced the multilevel Toeplitz and diagonal matrix-sequences, the only other class we need is that of zero--distributed matrix-sequences, whose definition depends on Definition \ref{def-distribution}.
	
	\begin{definition}\label{def:zero-distr}[Zero--distributed sequence]
		A matrix-sequence $\{Z_n\}_n$ such that $\{Z_n\}_n\sim_\sigma0$ is referred to as a zero-distributed sequence. In other words, $\{Z_n\}_n$ is zero-distributed if and only if
		\[ \lim_{n\to\infty}\frac1{n}\sum_{i=1}^nF(\sigma_i(Z_n))=F(0),\qquad\forall\,F\in C_c(\mathbb R). \]
	\end{definition}
In a different language, more common in the context of preconditioning and of the convergence analysis of (preconditioned) Krylov methods, a zero--distributed matrix-sequence is a sequence of matrices showing a (weak)
clustering at zero in the sense of the singular values (see e.g.\cite{MR3674485,TyZa-LAA} and references therein).

With the notaion indicating by $\|\cdot\|$ the spectral norm (i.e. the maximal singular value or equivalently the induced Euclidean norm) and by $\|\cdot\|_1$ the trace norm (i.e. the sum of all singular values), the following result holds true \cite{MR3674485}.	
	\begin{theorem}\label{th:glt}
	\begin{description}
		\item[{GLT\,1.}] If $\{A_n\}_n\sim_{\rm GLT}\kappa$ then $\{A_n\}_n\sim_\sigma\kappa$. If $\{A_n\}_n\sim_{\rm GLT}\kappa$ and the matrices $A_n$ are Hermitian then $\{A_n\}_n\sim_\lambda\kappa$.
		\item[{GLT\,2.}] If $\{A_n\}_n\sim_{\rm GLT}\kappa$ and $A_n=X_n+Y_n$, where
		\begin{itemize}
			\item every $X_n$ is Hermitian,
			\item $\|X_n\|,\,\|Y_n\|\le C$ for some constant $C$ independent of $n$,
			\item $n^{-1}\|Y_n\|_1\to0$,
		\end{itemize}
		then $\{A_n\}_n\sim_\lambda\kappa$.
		\item[{GLT\,3.}] We have
		\begin{itemize}
			\item $\{T_n(f)\}_n\sim_{\rm GLT}\kappa(x,\theta)=f(\theta)$ if $f\in \mathbb{L}^1([-\pi,\pi]^d)$,
			\item $\{D_n(a)\}_n\sim_{\rm GLT}\kappa(x,\theta)=a(x)$ if $a:[0,1]^d\to\mathbb C$ is Riemann-integrable, 
			\item $\{Z_n\}_n\sim_{\rm GLT}\kappa(x,\theta)=0$ if and only if $\{Z_n\}_n\sim_\sigma0$.
		\end{itemize}
		\item[{GLT\,4.}] If $\{A_n\}_n\sim_{\rm GLT}\kappa$ and $\{B_n\}_n\sim_{\rm GLT}\xi$ then
		\begin{itemize}
			\item $\{A_n^*\}_n\sim_{\rm GLT}\overline\kappa$,
			\item $\{\alpha A_n+\beta B_n\}_n\sim_{\rm GLT}\alpha\kappa+\beta\xi$ for all $\alpha,\beta\in\mathbb C$,
			\item $\{A_nB_n\}_n\sim_{\rm GLT}\kappa\xi$.
		\end{itemize}
		\item[{GLT\,5.}] If $\{A_n\}_n\sim_{\rm GLT}\kappa$ and $\kappa\ne0$ a.e.\ then $\{A_n^\dag\}_n\sim_{\rm GLT}\kappa^{-1}$.
	\end{description}
	\end{theorem}
A more general and more advanced result regarding item \textbf{GLT2} can be found in \cite{BaSe,GoSe}, even if for our purposes item \textbf{GLT2} is sufficient in our setting.

\subsection{Coco--Russo method in $2\mathbb{D}$: Dirichlet problem in a square domain}\label{ssec:method 2D}

We consider the following Dirichlet problem:
\begin{equation}
\begin{cases}
- u_{xx} - u_{yy} = f & \ \textrm{in} \ \Omega\\
u = g & \ \textrm{in} \ \partial \Omega,
\end{cases}
\label{2D_Problem}
\end{equation}
where $\Omega = [0,1] \times [a,1] \subset [0,1]^2$, $f,g: \Omega \rightarrow \mathbb{R}$ are assigned functions and $u: \Omega \rightarrow \mathbb{R}$ is the unknown function.

\begin{figure}[htbp]
\centering
\includegraphics[scale=.5]{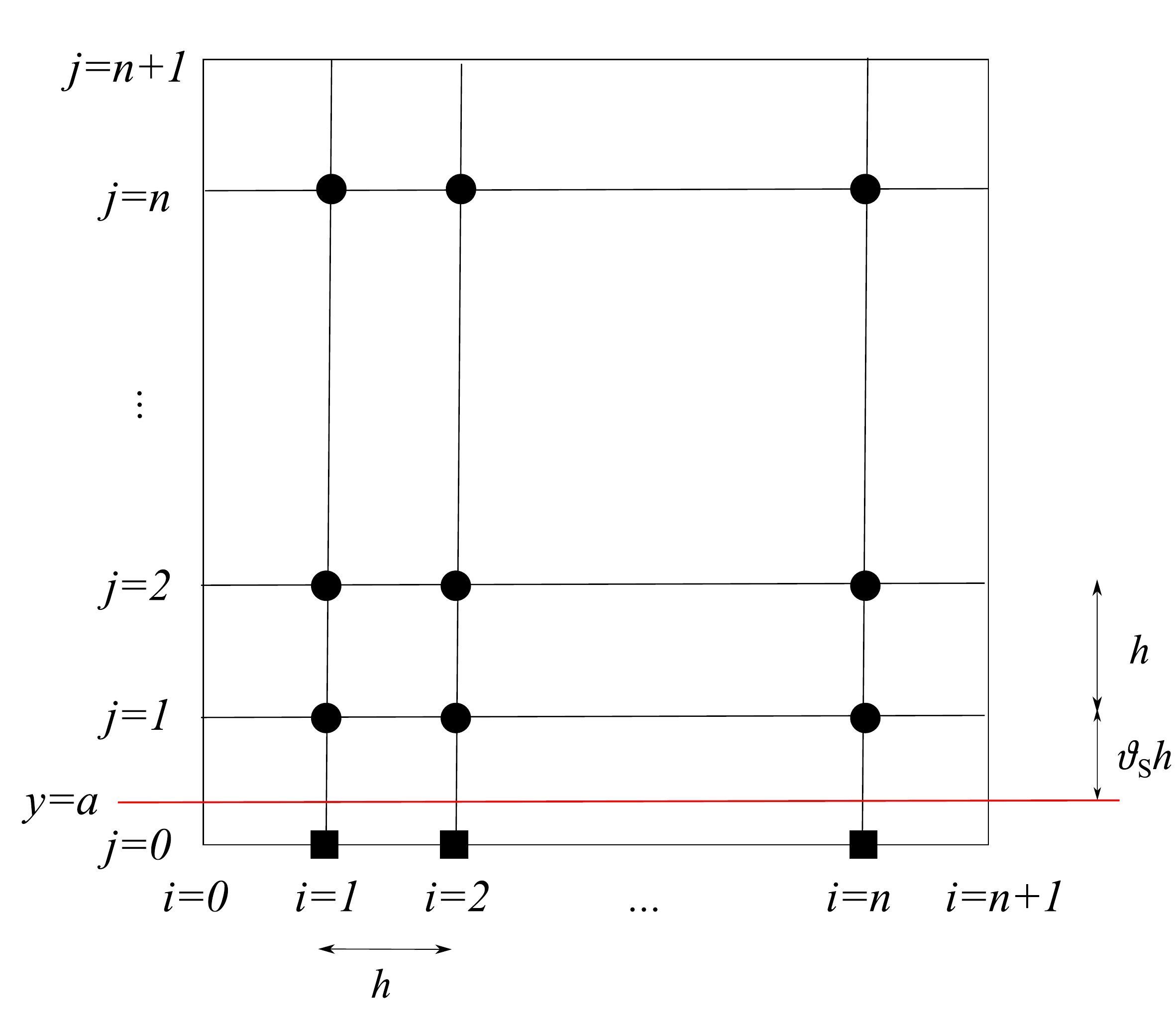}
\caption{Discretization of the $2\mathbb{D}$ domain. Full circles are the $n^2$ inside grid points, while full squares are the $n$ ghost points. Linear extrapolation is used to define the ghost values $u_{i,0}$ from $u_{i,1}$ and the boundary values $g(x_i,0)$, for $i=1,\ldots,n$.}
\label{domain2D}
\end{figure}

The square $[0,1]^2$ is discretized through a uniform Cartesian grid with $(n+2)^2$ grid points $(x_i,y_j) = (i\,h,j\,h)$, for $i,j=0,\ldots,n+1$, where $h=1/(n+1)$. As in the $1\mathbb{D}$ case, let $0<a<h$ and call $\vartheta_\textsc{S} = (x_1-a)/h$ (see Fig.~\ref{domain2D}). The subscript $S$ stands for {\em south}, since the boundary $y=a$ is the bottom side of the domain. A similar approach can be followed in the other cases.

The elliptic equation $-\Delta u =f$ of problem~\eqref{2D_Problem} is discretized by central finite difference on internal grid points, with eliminated boundary conditions on the boundaries $x=0$, $x=1$ and $y=1$. Then, for $2 \leq i \leq n-1$ and $1\leq j \leq n-1$ we have:
\[
\frac{4 u_{ij}-(u_{i-1\,j} + u_{i+1\,j} + u_{i\,j-1} + u_{i\,j+1})}{h^2} = f_{ij},
\]
while for  for $i=1$ and $j=1,\ldots,n-1$ we eliminate the boundary condition on $x=0$:
\[
\frac{4 u_{1\,j}-(u_{2\,j} + u_{1\,j-1} + u_{1\,j+1})}{h^2} = f_{ij} + \frac{g(0,y_j)}{h^2}.
\]
Similarly, we elimiante the boundary conditions on $x=1$ and $y=1$.
The boundary condition on $y=a$ is discretized by linear interpolation:
\[
\vartheta_\textsc{s} u_{i,0} + (1-\vartheta_\textsc{s}) u_{i,1} = g(x_i,a), \quad \text{ for } i=1,\ldots, n.
\]
Overall, there are $n^2$ inside grid points $(x_i,y_j)$ for $i,j=1,\ldots,n$ and $n$ ghost points for $(x_i,0)$ for $i=1\ldots,n$.

Using a total lexicographical order, the matrix of coefficients that we obtain is a $2$-level matrix with the following structure:
\begin{equation*}
A_h =
\left(
\begin{array}{cccccccccccccccc}
\vartheta_\textsc{s}\mathbb{I}_n & | & (1-\vartheta_\textsc{s})\mathbb{I}_n & | & & | &  & | &  & | & & \\
\hline
B & | & G & | & B & | &  & | &  & | & & \\
\hline
& | & B & | & G & | & B & | &  & | & & \\
\hline
& | &  & | & \ddots & | & \ddots & | & \ddots & | & & \\
\hline
& | &  & | &  & | & B & | & G & | & B &    \\
\hline
& | &  & | &  & | &  & | & B & | & G &     \\
\end{array}
\right),
\end{equation*}
where
\begin{equation*}
\vartheta_\textsc{s}\mathbb{I}_n \in \mathbb{R}^{n \times n},
\end{equation*}
\begin{equation*}
(1-\vartheta_\textsc{s})\mathbb{I}_n \in \mathbb{R}^{n \times n},
\end{equation*}
\begin{equation*}
B=-\frac{1}{h^2}\mathbb{I}_n \in \mathbb{R}^{n \times n},
\end{equation*}
and
\begin{equation*}
G =
\frac{1}{h^2}\begin{pmatrix}
4 & -1 &  &  &    \\
  -1 & 4 & -1 &  &    \\
   & \ddots& \ddots & \ddots &      \\
   &  & -1&4 & -1    \\
   &  &  & -1 & 4  \\
\end{pmatrix} \in \mathbb{R}^{n\times n},
\end{equation*}
$A_h$ has $n$ blocks of $G$, so $A_h\in\mathbb{R}^{n(n+1)\times n(n+1)}$.

\subsection{Spectral analysis in $1\mathbb{D}$ and in $2\mathbb{D}$}
	\label{ssec:spectral 2D}

Having in mind the notations of Subsection \ref{ssec:multilevel}, the matrix $A_h$ can be decomposed in the following way
\begin{align}
A_h&=\frac{1}{h^2}\left[T_{\mathbf n}(f) +
X_{\mathbf n}\right]
\end{align}
where $\mathbf n=(n+1,n)$, the size of $A_h$ is $N(\mathbf n)=n(n+1)$,
\begin{equation}\label{2-level T}
T_{\mathbf n}(f) = T_{n+1}(2-2\cos(\theta))\otimes\mathbb{I}_n+\mathbb{I}_{n+1}\otimes T_{n}(2-2\cos(\theta)),
\end{equation}
$T_{k}(2-2\cos(\theta))$ is a Toeplitz matrix, already used in the $1\mathbb{D}$ case in Section~\ref{sec:pb 1d}, and
\begin{align}\label{dec X_n}
X_{\mathbf n}&=\left[\begin{array}{c|c|c}
T_n\left(h^2\vartheta_S-4+2\cos(\theta)\right) &T_n\left(h^2(1-\vartheta_S)+1\right)&\mathbf{0}_{n\times n(n-1) }\\
\hline
\mathbf{0}_{n^2\times n}&\mathbf{0}_{n^2\times n}&\mathbf{0}_{n^2\times n(n-1)}
\end{array}\right].
\end{align}
Of course, taking into account relation (\ref{d-level T}) with $d=2$ and (\ref{2-level T}), the function $f$ is bivariate and can be written as
\[
f(\theta_1,\theta_2)=4-2\cos(\theta_1)-2\cos(\theta_2).
\]
Therefore by using item \textbf{GLT1} in Theorem \ref{th:glt} we have
\[
\{T_{\mathbf n}(f)\}\sim_{\rm GLT} f
\]
in the sense of of Subsection \ref{ssec:glt}, so that
\[
\{T_{\mathbf n}(f)\}\sim_{\sigma} f,
\]
according to Definition \ref{def-distribution}.
 Furthermore, since $T_{\mathbf n}(f)$ is Hermitian (in fact real symmetric) for any choice of the partial sizes, thanks to item \textbf{GLT1}, we deduce $\{T_{\mathbf n}(f)\}\sim_{\lambda} f$ as well.

Now, taking into account Definition \ref{def:zero-distr}, it is easy to see that $\{X_{\mathbf n}\}$ ia a zero--distributed matrix-sequence, sinply because its rank is bounded by $n$ and hence the number of nonzero singular values is at most $n=o(n(n+1))$ with $N(\mathbf n)=n(n+1)$ being the sinze of $X_{\mathbf n}$. Therefore by item \textbf{GLT3}
\[
\{X_{\mathbf n}\}\sim_{\rm GLT} 0,
\]
so that $\{h^2 A_h\}\sim_{\rm GLT} f$ by item \textbf{GLT4}, since both $\{T_{\mathbf n}(f)\}, \{X_{\mathbf n}\}$ are GLT matrix-sequences and $h^2 A_h =T_{\mathbf n}(f) +
X_{\mathbf n}$ for any choice of the partial sizes. Then, again by item \textbf{GLT1} we deduce
\[
\{h^2A_h\}\sim_{\sigma} f.
\]
However, $X_{\mathbf n}$ is non-Hermitian and therefore we cannot apply item \textbf{GLT1} for concluding $\{h^2A_h\}\sim_{\lambda} f$. However, this can be done by using item \textbf{GLT2}, as proven in the following lines both in $1\mathbb{D}$ and in $2\mathbb{D}$.

\begin{theorem}\label{sp distr 1D 2D}
With the notations used so far in $1\mathbb{D}$ we have
\begin{equation}\label{sp distr 1D}
\{h^2A_h\}\sim_{\lambda} 2-2\cos(\theta),
\end{equation}
while in $2\mathbb{D}$ we have
\begin{equation}\label{sp distr 2D}
\{h^2A_h\}\sim_{\lambda} 4-2\cos(\theta_1)-2\cos(\theta_2).
\end{equation}
\end{theorem}
\textbf{Proof}
In $1\mathbb{D}$ we recall the identity
\[
h^2A_h=T_n(2-2\cos(\theta))+ \mathbf{e}_1 \mathbf{v}_h^{\mathrm{T}}.
\]
Since $\mathbf{e}_1 \mathbf{v}_h^{\mathrm{T}}$ is a rank one matrix, it has a unique nozero singular value so that
\[
\|\mathbf{e}_1 \mathbf{v}_h^{\mathrm{T}}\|_1=\|\mathbf{e}_1 \mathbf{v}_h^{\mathrm{T}}\|=\|\mathbf{v}_h^{\mathrm{T}}\|_2
\]
and hence a trivial computation shows that
\[
\lim_{{n}\rightarrow \infty} \frac{\|\mathbf{e}_1 \mathbf{v}_h^{\mathrm{T}}\|_1}{n}=0.
\]
Therefore, by item \textbf{GLT2}, we infer that both the GLT matrix sequences $\{h^2 A_h\}, \{T_{n}(2-2\cos(\theta))\}$ share the same eigenvalue distribution function $2-2\cos(\theta)$, which is the GLT symbol, so that (\ref{sp distr 1D}) is proven.

In $2\mathbb{D}$, according to the $2$-level notation, we remind that
\[
h^2A_h=T_{\mathbf n}(4-2\cos(\theta_1)-2\cos(\theta_2)) +
X_{\mathbf n}, \ \ \ {\mathbf n}=(n+1,n).
\]
Now in the light of (\ref{dec X_n}) we deduce that
\[
\|X_{\mathbf n}\|_1\le \|T_n\left(h^2\vartheta_S-4+2\cos(\theta)\right)\|_1 +\|T_n\left(h^2(1-\vartheta_S)+1\right)\|_1.
\]
Now, using the fact that $2\pi\|T_n(g)\|_1\le n\int_{[-\pi,\pi}]|g(\theta)|\ d\theta$ (see \cite{se-ti}), we obtain
\[
\|X_{\mathbf n}\|_1\le {n}{2\pi}\int_{[-\pi,\pi}]|h^2\vartheta_S-4+2\cos(\theta)|\ d\theta + n (h^2(1-\vartheta_S)+1)
\]
and, as in the $1\mathbb{D}$ setting, if we divide by the size of $X_{\mathbf n}$ i.e. $n(n+1)$ we find
\[
\lim_{{\mathbf n}\rightarrow \infty}\frac{\|X_{\mathbf n}\|_1}{n(n+1)}=0.
\]

Consequently, again by item \textbf{GLT2}, we deduce that both the GLT matrix sequences $\{h^2 A_h\}, \{T_{\mathbf n}(4-2\cos(\theta_1)-2\cos(\theta_2))\}$ share the same eigenvalue distribution function $4-2\cos(\theta_1)-2\cos(\theta_2)$, which is the GLT symbol, and hence (\ref{sp distr 2D}) is proven.  \hfill $\bullet$
\ \\

The previous result shows a spectral distribution as nonnegative functions both in $1\mathbb{D}$ and $2\mathbb{D}$. More precisely, looking at the  range of the spectral symbols, we deduce that $[0,4]$ is a cluster for the eigenvalues of $\{h^2 A_h\}$ in $1\mathbb{D}$, while $[0,8]$ is a cluster for the eigenvalues of $\{h^2 A_h\}$ in $2\mathbb{D}$.

This is nontrivial (and somehow unexpected), given the fact that the related corrections are non-Hermitian and possess only strictly negative eigenvalues and zero eigenvalues.

\section{Conclusions}\label{sec:final}

We have provided spectral and norm estimates for matrix sequences arising from the approximation of the Laplacian via the Coco--Russo method and we have validated them with a few numerical experiments. The analysis has involved several tools from matrix theory and in particular from the setting of Toeplitz operators and Generalized Locally Toeplitz matrix sequences. Open problems remain involving variable coefficients and non square domains: both cases can be handled form a spectral view point using the GLT machinery. In particular  when considering variable coefficients, the use of the diagonal sampling matrix-sequences allows to remain in GLT $*$-algebra, while the case of non square domains can be treated using the reduced GLT theory (see page 398-399 in \cite{GLT1} and Subsection 3.1.4 in \cite{GLT2}).

More involved is the case of the norm estimates of the inverse even in the case of a square in $2\mathbb{D}$. Below we present an idea in this direction.

Actually the decomposition (\ref{dec X_n}) suggests, as in the $1\mathbb{D}$ setting, the use of the Sherman--Morrison--Woodbury formula: we can set $A=T_{\mathbf n}(f)$, $X_{\mathbf n}= UCV$, $\mathbf n=(n+1,n)$, so that
\begin{align} \nonumber
U&=\begin{bmatrix}
\mathbb{I}_n\\
\hline
\mathbf{0}_{n^2\times n}
\end{bmatrix}\in\mathbb{R}^{n(n+1)\times n}\\ \nonumber
C&=\mathbb{I}_n\in\mathbb{R}^{n\times n}\\ \nonumber
V&=\left[T_n\left(h^2\vartheta_S-4+2\cos(\theta)\right) | T_n\left(h^2(1-\vartheta_S)+1\right) | \mathbf{0}_{n\times n(n-1)}\right]\in\mathbb{R}^{n\times n(n+1)}\\ \nonumber
&=[V_1| V_2|\mathbf{0}_{n\times n(n-1)}].
\end{align}
Hence
\begin{align} \nonumber
(A+UCV)^{-1}=A^{-1}-A^{-1}U(C^{-1}+VA^{-1}U)^{-1}VA^{-1}
\end{align}
and thus $A_h^{-1}=h^2(A+UCV)^{-1}=h^2\left(A^{-1}-A^{-1}U(C^{-1}+VA^{-1}U)^{-1}VA^{-1}\right)$, with $C^{-1}=C=\mathbb{I}_n$

The previous reasoning can be useful and promising, since the entries of the inverse of $A=T_{\mathbf n}(f)$, $f(\theta_1,\theta_2)=4-2\cos(\theta_1)-2\cos(\theta_2)$, are explicitly known (see \cite{meurant}). However technical difficulties remain due to the complicate expression of the entries of $T_{\mathbf n}^{-1}(f)$: this task will be the subject of future investigations.

\section*{Acknowledgments}
Giovanni Russo and Stefano Serra-Capizzano are grateful to GNCS-INdAM for the support in the present research.
Giovanni Russo acknowledges support from the Italian Ministry of Instruction, University and Research (MIUR), PRIN Project 2017 (No. 2017KKJP4X entitled Innovative numerical methods for evolutionary partial differential equations and applications).

\addcontentsline{toc}{chapter}{References}
\bibliographystyle{abbrv}
\bibliography{bibliography}

\end{document}